\documentclass[graybox]{svmult}

\usepackage{type1cm}        
\usepackage{makeidx}         
\usepackage{graphicx}        
\usepackage{multicol}        
\usepackage[bottom]{footmisc}

\usepackage{newtxtext}       %
\usepackage[varvw]{newtxmath}       

\makeindex             

\begin{document}

\title*{Intrinsic polynomial squeezing for Balakrishnan-Taylor beam models}
\author{E. H. Gomes Tavares, M. A. Jorge Silva, V. Narciso and A. Vicente}
\institute{E. H. Gomes Tavares \at State University of Londrina, 86057-970, Londrina, PR, Brazil,\\ \email{eduardogomes7107@gmail.com}
\and M. A. Jorge Silva \at  State University of Londrina, 86057-970, Londrina, PR, Brazil. \\
\email{marcioajs@uel.br}
\and V. Narciso \at  State University of Mato Grosso do Sul, 79804-970, Dourados, MS, Brazil.\\ \email{vnarciso@uems.br}
\and A. Vicente \at  Western Paraná State University, 85819-110, Cascavel, PR, Brazil. \\
\email{andre.vicente@unioeste.br}}
%
%
\maketitle

\vskip-2.0cm


\abstract{We explore the  energy decay properties	related to a model in extensible beams with the so-called  {\it energy damping}. 
We investigate the influence of the nonloncal damping coefficient in the stability of the model. We prove, for the first time, that the corresponding energy functional is squeezed by polynomial-like functions involving the power of the  damping coefficient, which  arises intrinsically from the  Balakrishnan-Taylor beam models.
As a consequence, it is shown that  such models with nonlocal energy damping are never exponentially stable in its essence.}

\section{Introduction}
\label{sec-intro}
%

In 1989  Balakrishnan and Taylor \cite{BalakrishnanTaylor} derived some prototypes of  
vibrating extensible beams with the so-called  {\it energy damping}. Accordingly, the following  one dimensional beam equation is proposed
 \begin{eqnarray}\label{balak-taylor-eq1}
	\partial_{tt} u-2\zeta\sqrt{\lambda} \partial_{xx}u
	+ \lambda \partial_{xxxx} u
	-\alpha\left[\int_{-L}^{L}\big(\lambda |\partial_{xx} u|^2 + |\partial_t u|^2\big)dx\right]^q    \partial_{xxt} u = 0,
\end{eqnarray}
where $u=u(x,t)$ represents the transversal deflection of a beam with length  $2L>0$ in the rest position, $\alpha>0$ is a damping coefficient, $\zeta$ is a constant appearing in  Krylov-Bogoliubov's approximation, $\lambda>0$ is  related to  mode frequency   and  spectral density of external forces, and $q=2(n+\beta)+1$ with $n\in\mathbb{N}$ and $0\leq \beta < \frac{1}{2}$.  We still refer to \cite[Sect. 4]{BalakrishnanTaylor} for several other beam equations taking into account nonlocal energy damping  coefficients, as well as  
\cite{Balakrishnan,Bass,Dowell,Hughes,MuMa, You,Zhang} for associated models. A normalized $n$-dimensional equation  corresponding to \eqref{balak-taylor-eq1} can be seen as follows
\begin{equation}\label{ncontext}
	\partial_{tt} u-\kappa \Delta u + \Delta^2 u- \alpha\left[\int_{\Omega}\left(|\Delta u|^2+|\partial_t u|^2\right)dx \right]^q\Delta \partial_{t}u=0,
\end{equation}
  where we denote $\lambda=1$ and $\kappa=2\zeta$; $\Omega$ may represent an open bounded of $\mathbb{R}^n$; and the symbols  $\Delta$ and $\Delta^2$ stand for the usual Laplacian and Bi-harmonic operators, respectively. Additionally, in order to see the problem within the frictional context of dampers,
  we rely on
  materials whose viscosity can be essentially seen as friction between moving solids. In this way, besides
  reflecting on a more challenging model (at least) from the stability point of view, one may metaphysically
  supersede  the viscous damping in \eqref{ncontext} by a nonlocal frictional one so that we  cast the model
 \begin{equation}\label{ncontext-weak}
 	\partial_{tt} u-\kappa \Delta u + \Delta^2 u + \alpha\left[\int_{\Omega}\left(|\Delta u|^2+|\partial_t u|^2\right)dx\right]^q  \partial_{t}u=0.
 \end{equation}  
   
The main goal of this paper is to explore the influence of the nonloncal damping coefficient in the stability of problem \eqref{ncontext-weak}. Unlike the existing literature on extensible beams with full viscous or frictional damping, we are going to see for the first time that the feature of the {\it energy} damping coefficient
\begin{equation}\label{energy-coef}
	 \mathcal{E}_q(t):=  \mathcal{E}_q(u,u_t)(t)= \left[\int_{\Omega}\left(|\Delta u(t)|^2+|\partial_{t} u(t)|^2\right)dx\right]^q, \quad  q>0, 
\end{equation}
not only prevents exponential decay, but also gives us a polynomial range in terms of $q$ whose energy is squeezed and goes to zero 
polynomially when time goes to infinity. More precisely, by noting that the corresponding energy functional is given by 
\begin{equation}\label{energy-full}
 {E}_{\kappa}(t):=  {E}_{\kappa}(u,u_t)(t)=  \int_{\Omega}\left(|\Delta u(t)|^2+|\partial_{t} u(t)|^2 + \kappa |\nabla u(t)|^2 \right)dx, \quad \kappa\geq0, 
\end{equation}
then it belongs to an  area of  variation between upper and lower polynomial limits as follows
\begin{equation}\label{polin-rate}
c_0 \, t^{-\frac{1}{q}} \lesssim 	{E}_{\kappa}(t) \lesssim C_0 \, t^{-\frac{1}{q}}, \quad t \to +\infty,
\end{equation}
for some constants $0<c_0\leq C_0$ depending on the initial energy ${E}_{\kappa}(0), \, \kappa\geq0$. Indeed, such a claim corresponds to  an intrinsic polynomial range of (uniform) stability and will follow as a consequence 
of a more general result that is  rigorous stated in Theorem \ref{Main-bounds}. See also Corollary \ref{cor-polin}. In particular, we can conclude that \eqref{ncontext-weak} is not exponentially stable   when dealing with weak initial data, that is, with solution in the standard energy space. See Corollary \ref{cor-exp-never}. 


In conclusion, Theorem \ref{Main-bounds} truly  reveals  the stability of the associated energy  ${E}_{\kappa}(t)$, which leads us to the concrete conclusions provided by Corollaries \ref{cor-polin}-\ref{cor-exp-never}, being pioneering results on the subject.
Due to technicalities in the well-posedness process, we shall work with $q\geq1/2$. In Section \ref{sec-wp} we prepare all notations and initial results. Then,
 all precise details on the stability results shall be given in  Section \ref{sec-stability}.
%


\subsection{Previous literature, comparisons and highlights}
  

In what follows, we are going to highlight that our approach and results are different or else provide generalized results, besides keeping more physical consistency in working exactly with \eqref{energy-coef} instead of modified versions of it. Indeed, there are at least three mathematical ways of attacking the energy damping coefficient  \eqref{energy-coef} along the equation \eqref{ncontext-weak} (or \eqref{ncontext}), namely:
\begin{itemize}
	\item[1.] Keeping the potential energy in \eqref{energy-coef}, but neglecting the kinetic one; 
	
	\item[2.] Keeping the kinetic energy in \eqref{energy-coef}, but neglecting the potential one;
	
	\item[3.]  Keeping both potential and kinetic energies, but considering them under the action of a  strictly (or not) positive  function $M(\cdot)$  as a  non-degenerate (or possibility degenerate) damping coefficient.
\end{itemize}
 
In the first case, equation \eqref{ncontext-weak} becomes to
 \begin{equation}\label{ncontext-weak-pot}
	\partial_{tt} u-\kappa \Delta u + \Delta^2 u + \alpha\left[\int_{\Omega}|\Delta u|^2 dx\right]^q  \partial_{t}u=0 \ \ \mbox{ in } \ \ \Omega\times(0,\infty).
\end{equation}  
This is, for sure, the most challenging case once the damping coefficient becomes now to a real degenerate  coefficient. In  \cite[Theorem 3.1]{cavalcante2-marcio-vando}, working on a bounded domain $\Omega$ with clamped boundary condition, it is proved the following with  $q=1$ in \eqref{ncontext-weak-pot}: 
{\it    for every $R > 0$, there exist constants $C_R=C(R)>0$ and $\gamma_R=\gamma(R)>0$ depending on $R$ such that 
 \begin{equation}\label{exp}
	E_{\kappa}(t) \ \leq C_R \, E_{\kappa}(0) \, e^{-\gamma_R t }, \quad  t>0,
\end{equation}
only  holds for every regular solution $u$ of   \eqref{ncontext-weak} with   initial data $(u_0, u_1)$  satisfying}
\begin{equation}\label{bdd-int}
\|(u_0,u_1)\|_{(H^4(\Omega)\cap H^2_0(\Omega))\times H^2_0(\Omega)}\leq R.
\end{equation}
We stress that \eqref{exp} only represents a  {\it local stability result} since it holds on  every ball with radius $R>0$ in the strong topology  $(H^4(\Omega)\cap H^2_0(\Omega))\times H^2_0(\Omega),$ but they are not independent of the initial data. Moreover, as observed by the authors in \cite{cavalcante2-marcio-vando}, the 
drawback of \eqref{exp}-\eqref{bdd-int} is that it could not be proved    in the weak topology $H_0^2(\Omega) \times L^2(\Omega)$, even taking initial data uniformly bounded in $H_0^2(\Omega) \times L^2(\Omega).$
 Although we recognized that our results for \eqref{ncontext-weak} can not be  fairly compared to such a result, we do can conclude by means of the upper and lower polynomial bounds \eqref{polin-rate} that 
 the estimate \eqref{exp} will  never be reached for weak initial data given in  $H_0^2(\Omega) \times L^2(\Omega)$. Therefore, 
our results act as complementary 
conclusions to \cite{cavalcante2-marcio-vando}
by clarifying such drawback  raised therein, and yet giving a different point of view of  stability by means of \eqref{polin-rate} and its consequences concerning problem \eqref{ncontext-weak}.

In the second case, equation \eqref{ncontext-weak} falls into
\begin{equation}\label{ncontext-weak-kin}
	\partial_{tt} u-\kappa \Delta u + \Delta^2 u + \alpha\left[ \int_{\Omega}|\partial_{t} u|^2 dx\right]^q  \partial_{t}u=0 \ \ \mbox{ in } \ \ \Omega\times(0,\infty).
\end{equation}  
Unlike the first case, here we have an easier setting because the kinetic damping coefficient provides a kind of monotonous (polynomial) damping whose computations to achieve \eqref{polin-rate} remain unchanged (and with less calculations). This means that all results highlighted previously still hold  for this particular case. In addition, they clarify what is precisely the stability result related to problems addressed in  \cite{zhao1,zhao2}, which in turn represent particular models of abstract damping given by \cite[Section 8]{Aloiui-Ben-Haraux}. In other words, in terms of stability, our methodology provides a way to show the existence of absorbing sets with  polynomial rate (and not faster than polynomial rate depending on $q$)  when dealing with generalized problems relate to \eqref{ncontext-weak-kin}, subject that is not addressed in \cite{zhao1,zhao2}. 

Finally, in the third case let us see   equations \eqref{ncontext}-\eqref{ncontext-weak} as follows
 \begin{equation}\label{ncontext-weakM}
	\partial_{tt} u-\kappa \Delta u + \Delta^2 u +  M \left(\int_{\Omega}\left(|\Delta u|^2+|\partial_t u|^2\right)dx\right)  A\partial_{t}u=0 \  \mbox{ in } \  \Omega\times(0,\infty),
\end{equation} 
where operator $A$ represents the Laplacian operator $A=-\Delta$  or else the identity one $A=I$. 
Thus, here we clearly have two subcases, namely, when $M(\cdot)\geq0$ is a non-degenerate or possibly degenerate function. For instance, when $M(s)= \alpha s^q, \, s\geq0,$ and $A=-\Delta$, then we go back to problem \eqref{ncontext}. For this (degenerate) nonlocal  strong damping situation with $q\geq1$, 
it is considered in \cite[Theorem 3.1]{jorge-narciso-vicente}  
an upper polynomial stability for the corresponding energy, which also involves a standard nonlinear source term. Nonetheless, we call the attention to the following prediction result provided in \cite[Theorem 4.1]{jorge-narciso-vicente} for  \eqref{ncontext} addressed on a 
bounded domain $\Omega$ with clamped boundary condition and  $q\geq1$: {\it By taking  finite initial energy $0<E_\kappa(0)<\infty$, then 
 $E_\kappa(t)$ given in  \eqref{energy-full} satisfies
 	\begin{equation}\label{false-exp}
 		E_\kappa(t)\le 3E_\kappa(0) e ^{-\delta \int_{0}^{t} \| 
 		u(s)\|^{2q}ds }, \quad   t> 0,
 \end{equation}
 where   $\delta=\delta(\frac{1}{E_\kappa(0)})>0$ is a constant   proportional to $1/E_\kappa(0)$. 	
}

Although the estimate \eqref{false-exp} provides a new result with an {\it exponential face}, it does not mean any kind of stability result. Indeed, it is only a peculiar estimate indicating that prevents exponential decay patterns as remarked in \cite[Section 4]{jorge-narciso-vicente}. In addition, it is worth pointing out that   our computations to reach the stability result for problem \eqref{ncontext-weak} can be easily adjusted to \eqref{ncontext}, even for $q\geq1/2$ thanks to a inequality provided in \cite[Lemma 2.2]{Aloiui-Ben-Haraux}. Therefore, through the polynomial range \eqref{polin-rate} we provide here a much more accurate stability result than the estimate expressed by   \eqref{false-exp}, by concluding indeed that both problems \eqref{ncontext} and \eqref{ncontext-weak} are never exponentially stable in the topology of the energy space.
 
On the other hand, in the non-degenerate case $M(s)>0, \, s\geq0$, but still taking  $A=-\Delta$, a generalized version of \eqref{ncontext-weakM} has been recently approached by \cite{sun-yang} in a context of {\it strong attractors}, that is, 
the existence of attractors
in the topology of more regular space than the weak phase space. In this occasion, the $C^1$-regularity for $M>0$ brings out the non-degeneracy of the damping coefficient, which in turn allowed them to reach interesting results on well-posedness,  regularity and long-time behavior of solutions over more regular spaces.  
Such assumption of positiveness for the damping coefficient has been also addressed by other authors for related problems, see e.g.  \cite{jorge-narciso-DIE,jorge-narciso-CDCS,jorge-narciso-EECT}. From our point of view, in spite of representing a nice case, the latter does not portray the current situation of this paper so that we do not provide more detailed comparisons with such a non-degenerate problems, but we refer to \cite{cavalcante2-marcio-vando,jorge-narciso-DIE,jorge-narciso-CDCS,jorge-narciso-EECT,jorge-narciso-vicente,sun-yang} for a nice survey on this kind of non-degenerate damping coefficients. Additionally, we note that the suitable case of  non-degenerate damping coefficient   $M(s)>0, \, s\geq0$, and    $A=I$ in \eqref{ncontext-weakM} has not been considered in the literature so far and shall be concerned in another work by the authors in the future.

At light of the above statements, one sees e.g. when  $M(s)= \alpha s^q, \, s\geq0,$ and $A=I$, then problem \eqref{ncontext-weakM} falls into \eqref{ncontext-weak},  
being a problem not yet addressed in the literature that brings out a new branch of studies for such a nonlocal (possibly degenerate) damped problems, and also justifies all new stability results previously specified.

\section{The problem and well-posedness}\label{sec-wp}

Let us consider again the beam model with energy damping
\begin{equation}\label{main}
	\partial_{tt}u+\Delta^2 u -\kappa\Delta u+\alpha\left[\int_{\Omega}\left(|\partial_tu|^2+|\Delta u|^2\right)\, dy\right]^q \partial_t u=0 \   \mbox{ in }   \ \Omega\times(0,\infty),
\end{equation}
with clamped boundary condition 
\begin{equation}\label{initial-boundary}
	u=\frac{\partial u}{\partial \nu}=0 \   \mbox{ on }   \  \partial \Omega \times  [0,\infty),
\end{equation}
and initial data
\begin{equation}\label{initial-data}
	u(x,0)=u_0(x),\quad \partial_t u(x,0)=u_1(x), \quad x \in \Omega.
\end{equation}

To address problem \eqref{main}-\eqref{initial-data}, we introduce the  Hilbert  phase space (still called {\it energy space})
$$
\mathcal{H}:=H^2_0(\Omega) \times L^2(\Omega), 
$$
equipped with the inner product
$
\left\langle z^1,z^2\right\rangle_{\mathcal{H}}:=\left\langle\Delta u^1,\Delta u^2\right\rangle +\left\langle v^1,v^2\right\rangle$ for $z^i=(u^i,v^i)\in \mathcal{H},\,\, i=1,2,$
and norm
$\|z\|_{\mathcal{H}}=\left(\|\Delta u\|^2+\|v\|^2 \right)^{1/2},$ for  $z=(u,v) \in \mathcal{H},$
where $\left\langle u,v\right\rangle:=\displaystyle \int_{\Omega}uv\,dx$, $\|u\|^2:=\left\langle u,u \right\rangle$ and $\|z\|^2_{\mathcal{H}}:=\left\langle z,z\right\rangle_{\mathcal{H}}$.
\smallskip

In order to stablish the well-posedness of \eqref{main}-\eqref{initial-data}, we define the vector-valued function $
 z(t):=(u(t),v(t))$, $t\geq0,$ with 
 $v=\partial_{t}u$.  Then we can rewrite system \eqref{main}-\eqref{initial-data} as the following first order abstract problem
 \begin{eqnarray}\label{abstrac-cauchy}
 \left\{\begin{array}{l}
 \partial_t z =\mathcal{A}z +\mathcal{M}(z), \quad t>0, \medskip \\
 z(0)=(u_0,u_1):=z_0,
 \end{array}\right.
 \end{eqnarray}
where   $\mathcal{A}:\mathcal{D}(\mathcal{A})\subset\mathcal{H}\to\mathcal{H}$ is the linear  operator given by
\begin{equation}\label{def_A}
\mathcal{A}z = ( v , - \Delta^2 u),
\quad    \mathcal{D}(\mathcal{A}):=H^4(\Omega)\cap  H^2_0(\Omega),
\end{equation}
and $\mathcal{M}:\mathcal{H}\to\mathcal{H}$ is the nonlinear operator
\begin{equation}\label{def_B}
\mathcal{M}(z)= (0, \kappa \Delta u-\alpha\|z\|^{2q}_{\mathcal{H}}v), \quad z=(u,v)\in\mathcal{H}.
\end{equation}
Therefore, the existence and uniqueness of solution to the system \eqref{main}-\eqref{initial-data} relies on the study of problem (\ref{abstrac-cauchy}). Accordingly, we have the following well-posedness result.

\color{black}

\begin{theorem} \label{theo-existence} 
Let $\kappa,\alpha \geq 0$ and $q\geq \frac{1}{2}$ be given constants.  If $z_0\in \mathcal{H}$, then  \eqref{abstrac-cauchy} has a unique mild   solution  $z$ in the class
$z\in C([0,\infty),\mathcal{H}).$

In addition, if $z_0 \in \mathcal{D}(\mathcal{A})$, then $z$  is a regular solution lying in the class  \begin{equation*}
	z\in C([0,\infty),\mathcal{D}(\mathcal{A}))\cap  C^1([0,\infty),\mathcal{H}).
\end{equation*} 
 
\end{theorem}
\begin{proof}
To show the local version of the first statement, it is enough to prove that $\mathcal{A}$  given in  (\ref{def_A})  is the infinitesimal generator of a $C_{0}$-semigroup of contractions $e^{\mathcal{A}t}$ (which is very standard) and $\mathcal{M}$ set in  (\ref{def_B}) is locally Lipschitz   on $\mathcal{H}$ which will be done next. Indeed, let $r>0$ and $z^1,z^2 \in \mathcal{H}$ such that $\max\{\|z^1\|_{\mathcal{H}},\|z^2\|_{\mathcal{H}}\} \leq r$. We note that  
\begin{equation}\label{Sem_I}
 \left\|\|z^1\|^{2q}_{\mathcal{H}}v^1-\|z^2\|^{2q}_{\mathcal{H}}v^2\right\|\leq 
\left[\|z^1\|^{2q}_{\mathcal{H}}+\|z^2\|^{2q}_{\mathcal{H}}\right]\|v^1-v^2\|+\left|\|z^1\|^{2q}_{\mathcal{H}}-\|z^2\|^{2q}_{\mathcal{H}}\right|\|v^1+v^2\|.
\end{equation}
The first term on the right side of \eqref{Sem_I} can be estimated by
\begin{eqnarray*}
	\left[\|z^1\|^{2q}_{\mathcal{H}}+\|z^2\|^{2q}_{\mathcal{H}}\right]\|v^1-v^2\|\leq 2r^{2q}
	\|z^1-z^2\|_{\mathcal{H}}.
\end{eqnarray*}
Now, from a suitable inequality provided in  \cite{Aloiui-Ben-Haraux}\footnote{
See \cite[Lemma 2.2]{Aloiui-Ben-Haraux}: {\it Let $X$ be a normed space with norm $\|\cdot\|_X$. Then, for any $s\geq 1$ we have
		\begin{equation}\label{Haraux-variant}
			\left|\|u\|_X^s-\|v\|_X^s\right| \leq s\max\{\|u\|_X,\|v\|_X\}^{s-1}\|u-v\|_X, \quad \forall \,\, u,v \in X.
		\end{equation} }
}
we estimate the second term as follows
\begin{equation*}
\left|\|z^1\|^{2q}_{\mathcal{H}}-\|z^2\|^{2q}_{\mathcal{H}}\right|\|v^1+v^2\|\leq4qr^{2q}\|z^1-z^2\|_{\mathcal{H}}.
\end{equation*}
Plugging the two last estimates in \eqref{Sem_I}, we obtain 
$$ \left\|\|z^1\|^{2q}_{\mathcal{H}}v^1-\|z^2\|^{2q}_{\mathcal{H}}v^2\right\|_{\mathcal{H}} \leq 2(2q+1)r^{2q}\|z^1-z^2\|_{\mathcal{H}}.$$
Thus,
$$\|\mathcal{M}(z^1)-\mathcal{M}(z^2)\|_{\mathcal{H}} \leq \left(\kappa+2(2q+1)\alpha r^{2q}\right)\|z^1-z^2\|_{\mathcal{H}},$$
and $\mathcal{M}$ is locally Lipschitz in $\mathcal{H}$.

Hence, according to Pazy \cite[Chapter 6]{Pazy}, if $z_0 \in \mathcal{H}$ ($z_0 \in D(\mathcal{A})$), there exists a time $t_{\max} \in (0,+\infty]$ such that \eqref{abstrac-cauchy} has a unique mild (regular) solution 
$$z \in C([0,t_{\max}),\mathcal{H})\,\, (z \in C([0,t_{\max}),D(\mathcal{A}))\cap  C^1([0,t_{\max}),\mathcal{H})).$$
Moreover, such time $t_{\max}$ satisfies either the conditions $t_{\max}=+\infty$ or else $t_{\max}<+\infty$ with
\begin{equation}\label{blow-up}
\lim_{t\to t_{\max}^-}\|z(t)\|_{\mathcal{H}}=+\infty.
\end{equation} 

In order to show that $t_{\max}=+\infty$, we consider $z_0 \in D(\mathcal{A})$ and the corresponding regular solution $z$ of \eqref{abstrac-cauchy}. Taking the inner product in $\mathcal{H}$ of \eqref{abstrac-cauchy} with $z$, we obtain
\begin{eqnarray}\label{energy_relation_AA}
\frac{1}{2}\frac{d}{dt}\left[\|z(t)\|^2_{\mathcal{H}}+\kappa\|\nabla u(t)\|^2\right]+\alpha  \|z(t)\|^{2q}_{\mathcal{H}}\|\partial_t u(t)\|^{2}=0 \quad t \in [0,t_{\max}).
\end{eqnarray}
Integrating (\ref{energy_relation_AA}) over $(0,t), \,t \in [0,t_{\max})$, we get
$$
\|z(t)\|_{\mathcal{H}}\leq (1+c'\kappa)^{1/2}\|z_0\|_{\mathcal{H}}, \quad  t \in [0,t_{\max}).
$$
Here, the constant $c'>0$ comes from the embedding $H^2_0(\Omega) \hookrightarrow H_0^1(\Omega)$. The last estimate contradicts \eqref{blow-up}. Hence, $t_{max}=+\infty$. Using a limit process, one can conclude the same result for mild solutions.

The proof of Theorem \ref{theo-existence} is then complete.
\end{proof}

\section{Lower-upper polynomial energy's bounds}\label{sec-stability}
 
By means of the notations introduced in Section \ref{sec-wp}, we recall that the energy functional corresponding to problem \eqref{main}-\eqref{initial-data} can be expressed by
\begin{equation}\label{energy}
E_{\kappa}(t)=\frac{1}{2}\left[\|(u(t),\partial_t u(t))\|^2_{\mathcal{H}}+\kappa\|\nabla u(t)\|^2\right], \quad t\geq 0.
\end{equation}

Our main stability result reveals that $E_{\kappa}(t)$ is squeezed  by decreasing polynomial functions as follows.

 \begin{theorem}\label{Main-bounds}
Under the assumptions of Theorem  $\ref{theo-existence} $,   
there exists an increasing function $\mathcal{J}:\mathbb{R}^+\to \mathbb{R}^+$ such that the energy $E_{\kappa}(t)$ satisfies
\begin{eqnarray}\label{inequality_pricipal}
 	\left[\,2^{q+1}\alpha q t+\big[E_{\kappa}(0)\big]^{-q}\,\right]^{-1/q}\le E_{\kappa}(t)\le\left[\frac{q}{\mathcal{J}\left(E_{\kappa}(0)\right)}(t-1)^{+}+\big[E_{\kappa}(0)\big]^{-q}\right]^{-1/q},
 	\end{eqnarray}
 	for all $t>0$, where we use the standard notation $s^+:=(s+|s|)/2$. 
 \end{theorem} 
\begin{proof}
Taking the scalar product in $L^2(\Omega)$ of \eqref{main} with $\partial_t u$, we obtain
\begin{eqnarray}\label{absorbing_CC}
\frac{d}{dt}E_{\kappa}(t)= -\alpha||(u(t),\partial_t u(t))||^{2q}_{\mathcal{H}}\|\partial_t u(t)\|^2, \quad t>0.
\end{eqnarray}

Let us prove the lower and upper estimates in (\ref{inequality_pricipal}) in the sequel.

\noindent {\it Lower bound.} We first note that 
$$
||(u(t),\partial_t u(t))||^{2q}_{\mathcal{H}}\|\partial_t u(t)\|^2 \leq  2^{q+1}\left[\,E_{\kappa}(t)\,\right]^{q+1},
$$
and replacing it in  (\ref{absorbing_CC}), we get
\begin{eqnarray}\label{non}
	\frac{d}{dt}E_{\kappa}(t) \geq -2^{q+1}\alpha\left[\,E_{\kappa}(t)\,\right]^{q+1}, \quad t>0.
\end{eqnarray}
Thus, integrating \eqref{non} and proceeding a straightforward computation, we reach the first inequality in (\ref{inequality_pricipal}).

\smallskip 
\noindent {\it Upper bound.}  Now, we are going to prove the second inequality of \eqref{inequality_pricipal}. To do so, we provide some proper estimates and then apply a Nakao's result (cf. \cite{Nakao0,Nakao}).

We start by noting that 
\begin{eqnarray}\label{absorbing_BB}
||(u(t),\partial_t u(t)||^{2q}_{\mathcal{H}}\|\partial_t u(t)\|^2\geq \|\partial_t u(t)\|^{2(q+1)},
\end{eqnarray}
and replacing  (\ref{absorbing_BB}) in (\ref{absorbing_CC}), we get
\begin{eqnarray}\label{absorbing_C'}
\frac{d}{dt}E_{\kappa}(t)+\alpha\|\partial_t u(t)\|^{2(q+1)}\leq 0, \quad t>0,
\end{eqnarray}
which implies that $E_{\kappa}(t)$ is non-increasing with $E_{\kappa}(t)\leq E_{\kappa}(0)$ for every $t>0$. Also, integrating (\ref{absorbing_C'}) from $t$ to $t+1$, we obtain
\begin{eqnarray}\label{C''}
\alpha\int_t^{t+1}\|\partial_t u(s)\|^{2(q+1)}\,ds\leq E_{\kappa}(t)-E_{\kappa}(t+1):= [\,D(t)\,]^2.
\end{eqnarray}
Using H\"older's inequality with $\frac{q}{q+1}+\frac{1}{q+1}=1$ and (\ref{C''}), we infer
\begin{eqnarray}
\int_t^{t+1}\|\partial_t u(s)\|^2ds\le \frac{1}{\alpha^{\frac{1}{q+1}}}[\,D(t)\,]^{\frac{2}{q+1}}.\label{G}
\end{eqnarray}
From the Mean Value Theorem for integrals, there exist $t_1\in [t,t+\frac{1}{4}]$ and $t_2\in [t+\frac{3}{4},t+1]$ such that
\begin{eqnarray}\label{eq37}
\|\partial_t u(t_i)\|^2\le 4\int_t^{t+1}\|\partial_t u(s)\|^2ds\le \frac{4}{\alpha^{\frac{1}{q+1}}}[\,D(t)\,]^{\frac{2}{q+1}},\,\, i=1,2.
\end{eqnarray}

On the other hand, taking 
 the scalar product in $L^2(\Omega)$ of \eqref{main} with $u$ and integrating the result over $[t_1,t_2]$, we have
\begin{eqnarray}\label{F}
\nonumber\int_{t_1}^{t_2}E_{\kappa}(s)\, ds & = &\int_{t_1}^{t_2}\|\partial_t u(s)\|^2\,ds+\frac{1}{2}\left[\left(\partial_t u(t_1),u(t_1)\right)-\left(\partial_t u(t_2),u(t_2)\right)\right]\\
&&-\frac{\alpha}{2}\int_{t_1}^{t_2}||(u(s),\partial_t u(s))||^{2q}_{\mathcal{H}}\left(\partial_t u(s),u(s)\right)ds .
\end{eqnarray}
Let us estimate the terms in the right side of \eqref{F}. Firstly, we note that through H\"older's inequality, 
(\ref{eq37}) and Young's inequality, we obtain  
\begin{eqnarray*}
\left|\left(\partial_t u(t_1),u(t_1)\right)-\left(\partial_t u(t_2),u(t_2)\right)\right|
&\leq & d\sum_{i=1}^2\|\partial_t u(t_i)\|\|\Delta u(t_i)\|\\
	&\leq & \frac{8 d}{\alpha^{\frac{1}{2(q+1)}}}[\,D(t)\,]^{\frac{1}{q+1}}
	\sup_{t_1\le s\le t_2}[E_{\kappa}(s)]^{1/2}\\
	&\le & 
	\frac{128\, d^2}{\alpha^{\frac{1}{q+1}}}
	[\,D(t)\,]^{\frac{2}{q+1}}
	+\frac{1}{8}\sup_{t_1\le s\le t_2}E_{\kappa}(s),
\end{eqnarray*}
where the constant $d>0$ comes from the embedding $H^2_0(\Omega) \hookrightarrow L^2(\Omega)$. Additionally, using that $E_{\kappa}(t)\leq E_{\kappa}(0)$, we have
$$
\| (u(t),\partial_t u(t))\|^{2q}_{\mathcal{H}}\leq  2^q\left[\,E_{\kappa}(t)\,\right]^q\leq 2^q\left[\,E_{\kappa}(0)\,\right]^q.
$$
From this and (\ref{G}) we also get
\begin{eqnarray*}
\left|\int_{t_1}^{t_2}||(u(s),\partial_t u(s))||^{2q}_{\mathcal{H}}\left(\partial_t u(s),u(s)\right)ds\right|
	&\leq &	\frac{2^{2q+3}d^2\left[\,E_{\kappa}(0)\,\right]^{2q}}{\alpha^{-\frac{q}{q+1}}}[\,D(t)\,]^{\frac{2}{q+1}}\\
&&+\frac{1}{8\alpha}	\sup_{t_1\le s\le t_2}E_{\kappa}(s).
\end{eqnarray*}

Regarding again (\ref{G}) and replacing the above estimates in (\ref{F}), we obtain
\begin{eqnarray}\label{A3}
\int_{t_1}^{t_2}E_{\kappa}(s)\,ds\leq\mathcal{K}\left(E_{\kappa}(0)\right)[\,D(t)\,]^{\frac{2}{q+1}}+\frac{1}{8}\sup_{t_1\le s\le t_2}E_{\kappa}(s),
\end{eqnarray}
where we set the function $\mathcal{K}$ as
$$\mathcal{K}(s):=\left[\frac{64\, d^2+1}{\alpha^{\frac{1}{q+1}}}+2^{(q+1)}d^2\alpha^{\frac{2q+1}{q+1}}s^{2q}\right]>0.$$

Using once more the Mean Value Theorem for integrals and the fact that $E_{\kappa}(t)$ is non-increasing, there exists $\zeta\in[t_1,t_2]$ such that
$$
\int_{t_1}^{t_2}E_{\kappa}(s)\,ds=E_{\kappa}(\zeta)(t_2-t_1)\ge \frac{1}{2}E_{\kappa}(t+1),
$$
and then 
$$
\sup_{t\le s\le t+1}E_{\kappa}(s)=E_{\kappa}(t)=E_{\kappa}(t+1)+[\,D(t)\,]^2 \leq 2\int_{t_1}^{t_2}E_{\kappa}(s)\,ds+[\,D(t)\,]^2.
$$
Thus, from this and (\ref{A3}), we arrive at
\begin{eqnarray*}
	\sup_{t\le s\le t+1}E_{\kappa}(s)&\le& [\,D(t)\,]^2+2\int_{t_1}^{t_2}E_{\kappa}(s)ds\\
	&\le&
	[\,D(t)\,]^2+2\mathcal{K}\left(E_{\kappa}(0)\right)[\,D(t)\,]^{\frac{2}{q+1}}+\frac{1}{4} \sup_{t\le s\le t+1}E_{\kappa}(s),
\end{eqnarray*}
and since $0<\frac{2}{q+1}\le 2$, we obtain
\begin{eqnarray}\label{xxxx}
\sup_{t\le s\le t+1}E_{\kappa}(s)\le\frac{4}{3}
[\,D(t)\,]^{\frac{2}{q+1}}\left[\,[\,D(t)\,]^{\frac{2q}{q+1}}+2\mathcal{K}\left(E_{\kappa}(0)\right)\right].
\end{eqnarray}
Observing that 
$	[\,D(t)\,]^{\frac{2q}{q+1}}\le \left[E_{\kappa}(t)+E_{\kappa}(t+1)\right]^{\frac{q}{q+1}}\le 2^{\frac{q}{q+1}} \left[ E_{\kappa}(0)\,\right]^{\frac{q}{q+1}},$
and denoting by 
\begin{equation}\label{Jei}
\mathcal{J}(s):= \left(\frac{4}{3}\right)^{q+1}\left[(2s)^{\frac{q}{q+1}}+2\mathcal{K}(s)\right]^{q+1}>0,	
\end{equation}
and also recalling the definition of $[D(t)]^2$ in (\ref{C''}), we obtain from (\ref{xxxx}) that 
$$
	\sup_{t\le s\le t+1}[E_{\kappa}(s)]^{q+1}\le
	\mathcal{J}\left(E_{\kappa}(0)\right)[\,E_{\kappa}(t)-E_{\kappa}(t+1)\,].\nonumber
$$
Hence, applying e.g.  Lemma 2.1 of  \cite{Nakao} with $E_{\kappa}=\phi$, $\mathcal{J}\left(E_{\kappa}(0)\right)=C_0,$ and $K=0$, we conclude
$E_{\kappa}(t)\le\left[\frac{q}{\mathcal{J}\left(E_{\kappa}(0)\right)}(t-1)^{+}+\frac{1}{\big[E_{\kappa}(0)\big]^{q}}\right]^{-1/q},$
which ends the proof of the second inequality in (\ref{inequality_pricipal}).

The proof of  Theorem \ref{Main-bounds} is therefore complete.
\end{proof} 
 
\begin{remark}{\rm 
	It is worth point out that we always have
\begin{equation}\label{ineq}
	\left[\,2^{2q+1}\alpha q t+\big[E_{\kappa}(0)\big]^{-q}\,\right]^{-1/q}\le \left[\frac{q}{\mathcal{J}\left(E_{\kappa}(0)\right)}(t-1)^{+}+\big[E_{\kappa}(0)\big]^{-q}\right]^{-1/q},
\end{equation}	
so that 	it makes sense to express $E_{\kappa}(t)$ between the inequalities in \eqref{inequality_pricipal}. Indeed, from the definition $\mathcal{J}$ in \eqref{Jei} one easily sees that 
$\mathcal{J}\left(E_{\kappa}(0)\right)\geq\frac{1}{2^{2q+1}\alpha},$
from where one concludes  \eqref{ineq} promptly.}
\end{remark}

\begin{corollary}\label{cor-polin}
	{\bf (Polynomial Range of Decay).}
Under the assumptions of Theorem $\ref{Main-bounds}$, 	
the energy functional $E_{\kappa}(t)$ defined in \eqref{energy} decays squeezed as follows
\begin{equation} 
	c_0 \, t^{-\frac{1}{q}} \lesssim 	{E}_{\kappa}(t) \lesssim C_0 \, t^{-\frac{1}{q}}  \ \ \mbox{ as } \ \ t \to +\infty,
\end{equation}
for some constants $0<c_0\leq C_0$ depending on the initial energy ${E}_{\kappa}(0).$ 

\smallskip

In other words, $E_{\kappa}(t)$ decays polynomially at rate $t^{-1/q}$ ($q\geq1/2$) as  $t \to +\infty$.
\qed
\end{corollary}

 \begin{corollary}\label{cor-exp-never}{\bf (Non-Exponential Stability).}
Under the assumptions of Theorem $\ref{Main-bounds}$, 
the energy   $E_{\kappa}(t)$ set in \eqref{energy} never decays exponentially as $e^{- at }$   ($a>0$)  as $t \to +\infty$. 
 	\qed
 \end{corollary}

\end{document}